\documentclass[12pt]{amsart}
\usepackage{amsfonts}
\usepackage{amscd}

\newtheorem{lemma}{Lemma}[section]

\newtheorem{remark}{Remark}[section]

\newcommand{\BN}{{\mathbb{N}}}
\newcommand{\BR}{{\mathbb{R}}}

\newcommand{\BT}{{\mathbb{T}}}

\newcommand{\gC}{\Gamma}
\newcommand{\gc}{\gamma}

\newcommand{\gS}{\Sigma}
\newcommand{\gO}{\Omega}

\newtheorem{prop}{Proposition}[section]
\newtheorem{thm}[prop]{Theorem}
\newtheorem{lem}[prop]{Lemma}
\newtheorem{cor}[prop]{Corollary}

\theoremstyle{definition}
\newtheorem{Ack}[prop]{Acknowledgments}

\newtheorem{defn}[prop]{Definition}
\newtheorem{rem}[prop]{Remark}
\newtheorem{exam}[prop]{Example}

\begin{document}
\author{E. Breuillard, T. Gelander}
\date{January 2002}

\title{On dense free subgroups of Lie groups      }

\begin{abstract}
We give a method for constructing dense and free subgroups in real Lie groups. In particular we show
that any dense subgroup of a connected semisimple real Lie group $G$
contains a free group on two generators which is still dense in $G$, and
that any finitely generated dense subgroup in a connected non-solvable Lie
group $H$ contains a dense free subgroup of rank $\leq 2\cdot \dim H$. As an application, we obtain a new and elementary proof of a conjecture of Connes and Sullivan on amenable
actions, which was first proved by Zimmer.
\end{abstract}

\maketitle

\section{Introduction}

The main purpose of this paper is to give a proof of the following statement:

\begin{thm}
\label{the1}\label{th1} Let $G$ be a connected semi-simple real Lie group
and $\Gamma $ a dense subgroup of $G$. Then $\Gamma $ contains two elements
which generate a dense free subgroup of $G$.
\end{thm}

We shall also give the following generalizations of this result:

\begin{thm}
\label{th2} Let $G$ be a connected real Lie group which is topologically
perfect (i.e. the commutator $[G,G]$ is dense in $G$). If its Lie algebra $%
\mathfrak{g}$ is generated by $l$ elements, then any dense subgroup contains
a dense free subgroup of rank $l$.
\end{thm}

\begin{thm}
\label{th3} Let $G$ be a connected non solvable real Lie group of dimension $%
d$. Then any finitely generated dense subgroup of $G$ contains a dense free
subgroup of rank $2d$.
\end{thm}

The finite generation assumption in \ref{th3} is crucial. Moreover, as it will be clear from the proofs, one can find a free
dense subgroup of rank $k$, for any integer $k\geq l$ in theorem \ref
{th2} (resp. $k\geq 2d$ in theorem \ref{th3}). For some
particular groups $G$ one can give a smaller bound for the rank of the dense
free subgroup than the bounds $l$ and $2d$ given in \ref{th2} and \ref
{th3}. We will give examples to
illustrate these facts. \newline

Theorem \ref{th1} was originally motivated by a conjecture in the theory of
amenable actions of groups on measure spaces. It answers a question first
raised by Carri\`{e}re and Ghys (see \cite{Ghys}) in a note from 1985 on
amenable actions. The authors learned about this question from A. Lubotzky. 
Combined with \cite{Ghys} theorem 2, it gives a direct
proof of the following theorem:

\begin{thm}
\label{C-S}Let $G$ be a real Lie group and $\Gamma $ a countable subgroup.
Then the action of $\Gamma $ on $G$ (or on $G/P$ for $P$ closed amenable) by left
translations is amenable if and only if the connected component of the
closure of $\Gamma $ in $G$ is solvable. 
\end{thm}

This result was first conjectured by Connes and Sullivan (see \cite{Ghys})
and subsequently proved by Zimmer in \cite{Zim} using the techniques of
super-rigidity theory. In their note, Carri\`{e}re and Ghys show that a
non-discrete free subgroup of a real Lie group cannot act amenably by left
translations on the Lie group (see \cite{Ghys}, theorem 2) and then give a
quick proof of theorem \ref{th1} in the case $G=SL_{2}(\Bbb{R})$, hence
proving theorem \ref{C-S} for $SL_{2}(\Bbb{R})$. \newline

A celebrated theorem of Tits (see \cite{T}), asserts that any
Zariski-dense subgroup of a semisimple algebraic group $G$ over a field of
characteristic zero contains a Zariski-dense free subgroup on two generators. Theorem \ref
{th1} can then be viewed as a kind of topological version of Tits'
theorem. The difficulty of the problem we are considering comes from the
fact that the free group obtained in Tits' proof is in general discrete.
In this paper, we give a practical method for constructing free generators of a free group 
in Zariski dense subsets of $G$ (see theorem \ref{main}). As an application, we then 
obtain theorems \ref{th1} to \ref{th3}. One can also easily derive from it the original 
statement of Tits' alternative. The proof relies on a careful study
of the so-called ''proximal'' elements in $\Gamma $ acting on a projective
space over some local field. Unlike in Tits' paper, where powers of suitable
semisimple elements are used to produce proximal elements, our method is
inspired from that of Abels-Margulis-Soifer (see \cite{abels}), which is
based on the Cartan decomposition of projective transformations.\newline

In section \ref{first}, we give a simple method for producing finitely
generated dense subgroups in connected real Lie groups. Section \ref{dyn} is
devoted to the study of proximal transformations in projective spaces over
local fields. We give some quantitative estimates that enable us to exhibit
proximal elements with nice properties. These two sections can be read independently. 
In section \ref{cons} we explain how
to find a suitable linear representation of $\gC$ and elements in $\Gamma $ which have a nice proximal action
on the corresponding projective space. Finally, the last section contains the proof of the three theorems above.%
\newline

Let us make one remark about terminology.

\begin{rem}
In this paper, we use both the usual real topology (or the Hausdorff
topology induced by a local field $k$) and the Zariski topology on the
groups considered. To avoid confusion, we shall add the prefix Zariski- to
any topological notion regarding the Zariski topology (Zariski-connected,
Zariski-dense, etc.). For the real topology, however, we shall plainly say
''dense'' or ''open'' without further notice. Note that the Zariski topology
on rational points does not depend on the field of definition, i.e. if $V$
is an algebraic variety defined over a field $K$ and if $L$ is any extension
of $K$, then the $K$-Zariski topology on $V(K)$ coincide with the trace of
the $L$-Zariski topology on it.
\end{rem}


\section{Generating dense subgroups in Lie groups}

\label{first}

We call a connected Lie group $G$ \textbf{topologically perfect} if its
commutator group $[G,G]$ is dense, or equivalently, if $G$ has no continuous
surjective homomorphism to the circle.

The following result shows that, in a topologically perfect group, elements
which lie near the identity generate a dense subgroup, unless they have some
algebraic reason not to. Similar statements were established by Kuranishi
(see \cite{kuranishi}).

\begin{thm}[Generating dense subgroups in topologically perfect groups]
\label{denseTP} Let $G$ be a connected topologically perfect real Lie group
with Lie algebra $\frak{g}$. Then there is an identity neighborhood $\Omega
\subset G$, on which $\log =\exp ^{-1}$ is a well defined diffeomorphism,
such that $g_{1},\dots ,g_{m}\in \Omega $ generate a dense subgroup whenever 
$\log (g_{1}),\dots ,\log (g_{m})$ generate $\frak{g}$.
\end{thm}

\begin{proof}
Recall that a Zassenhaus neighborhood of a real Lie group is an identity 
neighborhood $\gO$ for which the intersection $\gO\cap\gS$ with any discrete 
subgroup $\gS$ is contained in a connected nilpotent Lie subgroup. 
A classical theorem
of Kazhdan and Margulis (see \cite{Rag} 8.16) says that a 
Zassenhaus neighborhood always exists.
Checking the details of the proof of the Kazhdan-Margulis theorem, one can 
easily
verify that the identity neighborhood $\gO$ established there, satisfies the 
stronger property that its image under any surjective homomorphism $f$ is a
Zassenhaus neighborhood in the image group (this is because the properties
of the neighborhoods established in 8.18 and 8.20 in \cite{Rag} are inherited
to their images). Moreover, as in \cite{Rag}, one writes $\gO =\exp V$
for a suitable neighborhood $V$ of $0$ in $\frak{g}$, and then, a subset 
$\{ X_1,\ldots ,X_p\}\subset (df)(V)$ generates a nilpotent Lie subalgebra iff
$\exp (X_1),\ldots ,\exp (X_p)$ generate a nilpotent subgroup in $f(G)$.
We shall call such an $\gO$ a {\it strongly Zassenhaus neighborhood}.

Let $\gO$ be a strongly Zassenhaus neighborhood in $G$, and assume that it is
small enough so that $\log|_\gO$ is a well defined diffeomorphism. 
Pick $g_1,\dots ,g_m\in\gO$ for which $\log (g_1),\dots ,
\log (g_m)$ generate $\frak{g}$. Denote by $A$ the 
closure of $\langle g_1,\dots ,g_m\rangle$, and by $A^0$ its
identity component. We need to show that $A^0=G$. 

First observe that $A^0$ is normal in $G$, i.e. its Lie subalgebra 
$\mathfrak{a}$ is an ideal in $\mathfrak{g}$.
Indeed, since $\mathfrak{g}=\langle\log (g_i)\rangle$, it is enough to 
show that 
$\textrm{ad}\big(\log (g_i)\big) (\mathfrak{a})=\mathfrak{a}$ for any 
$1\leq i\leq m$, but
$$
 \textrm{ad}(\log g_i)(\mathfrak{a})=
 \big(\log \textrm{Ad}(g_i)\big) (\mathfrak{a}),
$$
and since $g_i\in A$, it normalizes $A^0$, i.e.
$\textrm{Ad}(g_i)(\mathfrak{a})=\mathfrak{a}$, which implies also that
$\big(\log \textrm{Ad}(g_i)\big) (\mathfrak{a})=\mathfrak{a}$.

Next consider the quotient map $f:G\to G/A^0$. 
Clearly, $f\big(\langle g_1,\dots ,g_m\rangle\big) =f(A)$ is discrete, and generated by
the set $\{ f(g_1),\dots ,f(g_m)\}$ which is contained in the 
Zassenhaus neighborhood $f(\gO )$. Therefore 
$f\big(\langle g_1,\dots ,g_m\rangle\big)$ is a nilpotent group, but this 
means that
$$
 \log f(g_i)=(df) (\log g_i),~i=1,\ldots ,m
$$ 
generate a nilpotent Lie 
algebra. As $\{ (df) (\log g_i)\}_{i=1}^m$ generate 
$\frak{g}/\frak{a}$, this implies that $G/A^0$ is nilpotent. 
Since $G$ is assumed to be topologically perfect, any nilpotent quotient is 
trivial. Thus $G=A^0$ and  $\langle g_1,\dots ,g_m\rangle$ is dense.
\end{proof}

In fact, it is enough to require that $\log (g_{1}),\dots ,\log (g_{m})$
generate a Lie subalgebra which corresponds to a dense Lie subgroup. This
phenomenon is illustrated by the following:

\begin{exam}
Let $\widetilde{\text{SL}}_{2}({\mathbb{R}})$ be the universal covering of $%
\text{SL}_{2}({\mathbb{R}})$, and let $a\in \widetilde{\text{SL}}_{2}({%
\mathbb{R}})$ be a central element of infinite order. Let $\alpha $ be an
irrational rotation in the circle group ${\mathbb{T}}$. Consider the group 
\begin{equation*}
G=\big(\widetilde{\text{SL}}_{2}({\mathbb{R}})\times {\mathbb{T}}\big)%
/\langle (a,\alpha )\rangle .
\end{equation*}
$G$ is an example of a topologically perfect non-perfect group. Its Lie
algebra $\frak{g}\cong \text{sl}_{2}({\mathbb{R}})\oplus {\mathbb{R}}$ is
not generated by $2$ elements, but the $2$-generated Lie subalgebra $\text{sl%
}_{2}({\mathbb{R}})$ corresponds to a dense Lie subgroup in $G$. Theorem \ref
{th2} guarantees only that any dense subgroup of $G$ contains a dense $F_{3}$%
, however, it is possible to show (by the same argument which proves theorem 
\ref{th1}) that there is always a dense $F_{2}$.
\end{exam}

However, in general, 2 generators are not enough. This phenomenon is
illustrated by the following:

\begin{exam}
Consider the complex Lie group $G=\text{SL}_{2}({\mathbb{C}})\cdot ({%
\mathbb{C}}^{2})^{5}$, with the ordinary action on $\text{SL}_{2}({\mathbb{C}%
})$ on each of the 5 copies of ${\mathbb{C}}^{2}$. It is easy to verify that

\begin{enumerate}
\item  $G$ is perfect, i.e. $[G,G]=G$.

\item  The exponential map is onto.

\item  Any 2 elements in $G$'s Lie algebra $\mathfrak{g}=\text{sl}_{2}({%
\mathbb{C}})\cdot ({\mathbb{C}}^{2})^{5}$ generate a Lie subalgebra for
which the closure of the corresponding Lie subgroup is proper.
\end{enumerate}

These properties imply that there is no 2-generated dense subgroup in $G$.
\end{exam}

It is well known that a real semisimple Lie algebra is generated by 2
elements (cf. \cite{kuranishi} or \cite{Bourb} VIII, 3, ex. 10). 
In fact $V=\{ (X,Y)\in\frak{g}\times\frak{g}:\langle X,Y\rangle \neq \frak{g} \}$ is a proper algebraic subset.
We therefore obtain:

\begin{thm}[Generating dense subgroups in semisimple groups]
\label{denseSS} Let $G$ be a connected semisimple Lie group. Then there
exists an identity neighborhood $\Omega \subset G$, and a proper exponential
algebraic subvariety $R\subset \Omega \times \Omega $ such that $\langle
x,y\rangle $ is dense in $G$ whenever $(x,y)\in \Omega \times \Omega
\setminus R$.
\end{thm}

Theorem \ref{denseTP} yields a lot of freedom in the procedure of generating
dense subgroups in topologically perfect groups. In fact, if $G$ is a
topologically perfect Lie group and $g_{1},\ldots ,g_{m}\in G$ are elements
close enough to $1$ for which $\log (g_{1}),\ldots ,\log (g_{m})$ generate
the Lie algebra $\frak{g}$, and if $U\subset G$ is a sufficiently small identity
neighborhood, then for any selections $h_{i}\in g_{i}\cdot U$, $\log
(h_{1}),\ldots ,\log (h_{m})$ generate $\frak{g}$, and $h_{1},\ldots ,h_{m}$
generate a dense subgroup in $G$. In particular, we obtain:

\begin{cor}
\label{01} If $G$ is a topologically perfect Lie group for which the Lie
algebra is generated by $m$ elements, then any dense subgroup in $G$
contains a $m$-generated dense subgroup.
\end{cor}

\begin{proof}
Let $g_1,\ldots ,g_m$ and $U$ be as in the last paragraph above. Clearly, if 
$\gC\leq G$ is a dense subgroup then $\gC\cap g_i\cdot U\neq\emptyset$. Pick
$h_i\in\gC\cap g_i\cdot U,~i=1,\ldots ,m$ then $\langle h_1,\ldots ,h_m\rangle$
is dense.
\end{proof}

The requirement that $G$ is topologically perfect is crucial for corollary 
\ref{01} as explained by the following remark:

\begin{rem}
If a connected Lie group $G$ is not topologically perfect then it has a
surjective homomorphism to the circle. Let $\Gamma \leq G$ be the pre-image
of the group of rational rotations. Clearly $\Gamma $ is dense in $G$ but,
as any finitely generated group of rational rotations is finite, $\Gamma $
has no finitely generated dense subgroup.
\end{rem}

This remark explains why, when $G$ is not topologically perfect, we shall
consider only \textit{finitely generated} dense subgroups. Then we can
control the number of generators, by the following reasoning.

\begin{prop}
\label{denseInFG} Let $G$ be a connected Lie group and $\Gamma \leq G$ a
finitely generated dense subgroup. Then $\Gamma $ contains a dense subgroup
on $2\dim G$ generators. If $G$ is compact then $\Gamma $ contains a dense
subgroup on $\dim G$ generators.
\end{prop}

\begin{proof}
Assume first that $G$ is abelian. Then $G=\BR^{d_1}\times\BT^{d_2}$ and we can
find $d_1$ elements in $\gC$ which generate a discrete cocompact subgroup.
Dividing by this subgroup, we may assume that $G$ is a torus. Then we 
argue by induction on $\dim G$. As $\gC$ is finitely generated and dense
it contains an element $\gc$ of infinite order. By replacing $\gc$ by some 
power $\gc^j$ if necessary, we obtain an element which generates a subgroup 
with connected closure $C$ of positive dimension. By induction, the proposition
holds for $G/C$. Lift arbitrarily to $G$ a set of $\dim G/C$ 
generators for a dense subgroup of the image of $\gC$ in $G/C$. Together
with $\gc$ they generate a dense subgroup in $G$.

For the general case, we argue as follows.
Consider the sequence $G_{i+1}=\overline{[G_i,G_i]}$, where $[G_i,G_i]$ is the commutator
group of $G_i$, starting at $G_0=G$. As $\dim G$ is finite this sequence
stabilizes at some finite step $k$. Then $G_k$ is a connected closed normal
topologically perfect subgroup. 
Moreover the commutator $[\gC,\gC]$ is a finitely generated dense subgroup in $G_1$, and 
similarly, the $i$'th commutator of $\gC$ is a finitely generated dense 
subgroup in $G_i$.   
Thus, by the above paragraph, we can find, for each $i$, $0\leq i<k$, a set $\gS_i$
of $2\dim G_i/G_{i+1}$ (or $\dim G_i/G_{i+1}$ in the compact case) elements in 
$G_i\cap\gC$ which projects to generators of a dense subgroup in the abelian
group $G_i/G_{i+1}$.
Additionally, by corollary \ref{01} there is a set $\gS_k$ of $\dim G_k$ 
elements in $\gC\cap G_k$ which generate a dense subgroup of $G_k$ (in the 
non-abelian compact case $G_k=G_1$ is semisimple and we can take $\gS_k$ of 
size 2). 

Clearly $\bigcup_{i=0}^k \gS_i$ generates a dense subgroup and has 
the required cardinality.
\end{proof}


\section{Projective transformations, proximality\label{dyn}, ping-pong}

Let $k$ be a local field equipped with an absolute value $|\cdot |$. We wish
to investigate the action of projective transformations in PSL$_{n}(k)$ on
the projective space $\Bbb{P}^{n-1}(k)$. We shall start by recalling the Cartan
decomposition, i.e. the $KAK$ decomposition in SL$_{n}(k)$, and introducing a
nice metric on the projective space $\Bbb{P}^{n-1}(k)$.

\begin{enumerate}
\item  Consider first the case when $k$ is archimedean, i.e. $k=\Bbb{R}$ or $%
\Bbb{C}$. We shall denote by $\left\| \cdot \right\| $ the canonical
euclidean (resp. hermitian) norm on $k^{n}$, i.e.$\left\| x\right\|
^{2}=\sum \left| x_{i}\right| ^{2}$ if $x=\sum x_{i}e_{i}$, where $%
(e_{1},...,e_{n})$ is the canonical basis of $k^{n}$. The Cartan
decomposition in SL$_{n}(k)$ reads 
\begin{equation*}
\text{SL}_{n}(k)=KAK
\end{equation*}
where, 
\begin{eqnarray*}
K &=&\text{SO}_{n}(\Bbb{R})\text{ or SU}_{n}(\Bbb{C})\text{~according to
whether~}k={\mathbb{R}}\text{~or~}{\mathbb{C}}, \\
\text{~and~}A &=&\{\text{diag}(a_{1},\ldots ,a_{n}):a_{1}\geq \ldots \geq
a_{n}>0,\prod a_{i}=1\}.
\end{eqnarray*}
Any element $g\in \text{SL}_{n}(k)$ can be decomposed as a product 
$g=k_{g}a_{g}k_{g}^{\prime }$, where $k_{g},k_{g}^{\prime }\in K$ and $%
a_{g}\in A$. We remark that $a_{g}$ is uniquely determined by $g$, but $%
k_{g},k_{g}^{\prime }$ are not. Even so, we shall use the subscript $g$ to
indicate the relation to $g$.\newline

\item  Next, consider the case where $k$ is non-archimedean with valuation
ring $\mathcal{O}_{k}$ and uniformizer $\pi $. We endow $k^{n}$ with the
canonical norm defined by $\left\| x\right\| =\max \left| x_{i}\right| ,$
where $x=\sum x_{i}e_{i}$ is the expansion of $x$ with respect to the
canonical basis $(e_{1},...,e_{n})$ of $k^{n}$. Then we get 
\begin{equation*}
\text{SL}_{n}(k)=KAK
\end{equation*}
\end{enumerate}

where 
\begin{eqnarray*}
K &=&\text{SL}_{n}(\mathcal{O}_{k}), \\
\text{~and~}A &=&\{\text{diag}(\pi ^{j_{1}},\ldots ,\pi ^{j_{n}}):j_{i}\in {%
\mathbb{Z}},j_{i}\leq  j_{i+1},\sum j_{i}=0\}.
\end{eqnarray*}
Any $g\in $SL$_{n}(k)$ can be decomposed as a product $g=k_{g}a_{g}k_{g}^{%
\prime }$, where $k_{g},k_{g}^{\prime }\in K$ and $a_{g}\in A$ (see \cite{PR}
page 150, or \cite{BT} 4.4.3.). Again $a_{g}$ is uniquely determined, but $%
k_{g},k_{g}^{\prime }$ are not.\newline

In both cases, the canonical norm on $k^{n}$ gives rise to the associated
canonical norm on $\bigwedge^{2}k^{n}$. We shall then define the \textit{%
standard metric} on $\Bbb{P}^{n-1}(k)$ by the formula 
\begin{equation*}
d\big( [v],[w]\big) =\frac{\left\| v\land w\right\| } {\left\|
v\right\|\cdot\left\| w\right\| }
\end{equation*}
This is well defined and satisfies the following properties:

\begin{itemize}
\item  $d$ is a distance on $\Bbb{P}^{n-1}(k)$ which induces the canonical
topology inherited from the local field $k$.

\item  $d$ is a ultra-metric distance if $k$ is non-archimedean, i.e. 
\begin{equation*}
d\big( [v],[w]\big)\leq \max \{d\big( [v],[u]\big) ,d\big( [u],[w]\big)\}
\end{equation*}
for any non-zero vectors $u,v$ and $w$ in $k^{n}$.

\item  If $f$ is a linear form $k^{n}\to k$, then for any non-zero
vector $v\in k^{n},$%
\begin{equation}
d\big( [v],[\ker f]\big) =\frac{\left| f(v)\right| }{\left\| f\right\| \cdot
\left\| v\right\| }  \label{lin}
\end{equation}

\item  the compact group $K$ acts by isometries on $\big(\Bbb{P}^{n-1}(k),d%
\big).$
\end{itemize}

In the sequel, we shall denote by $a_{1}(g),...,a_{n}(g)$ the coefficients of
the diagonal matrix $a_{g}$ corresponding to $g$ in the Cartan
decomposition. For further use, we note that in the above notations, for any
matrix $g\in $SL$_{n}(k)$, 
\begin{eqnarray*}
\left\| g\right\| &=&|a_{1}(g)| \\
\left\| \bigwedge\nolimits^{2}g\right\| &=&|a_{1}(g)a_{2}(g)|.
\end{eqnarray*}
\newline

Let us observe a few properties of projective transformations. For $g\in 
\text{SL}_{n}(k)$ we denote by $[g]$ the corresponding projective
transformation $[g]\in \text{PSL}_{n}(k)$.

\begin{lem}
\label{biLip} Every projective transformation is bi-Lipschitz on the entire
projective space for some constant depending on the transformation.
\end{lem}

\proof%
As the compact group $K$ acts by isometries on $\Bbb{P}^{n-1}(k)$, the $KAK$
decomposition allows us to assume that $g=a_{g}=diag(a_{1},...,a_{n})\in A$.
Thus, one only needs to check the easy fact that $a_{g}$ is $\left| \frac{%
a_{1}}{a_{n}}\right| ^{2}$-Lipschitz.%
\endproof%

\begin{defn}
Let $\epsilon \in (0,1)$. A projective transformation $[g]\in \text{PSL}%
_{n}(k)$ is called \textbf{$\epsilon $-contracting} if there exist a point $%
v_{g}\in {\mathbb{P}}^{n-1}(k),$ called an \textbf{attracting point} of $[g],
$ and a projective hyperplane $H_{g}$, called a \textbf{repulsive hyperplane}
of $[g]$, such that $[g]$ maps the complement of the $\epsilon $%
-neighborhood of $H_{g}\subset \Bbb{P}^{n-1}(k)$ into the $\epsilon $-ball
around $v_{g}$. We say that $[g]$ is $\epsilon $\textbf{-very contracting}
if both $[g]$ and $[g^{-1}]$ are $\epsilon $-contracting.
\end{defn}

Note that in general $v_{g}$ and $H_{g}$ are not necessarily unique.
Nevertheless, all the statements we shall make about them will be valid for
any choice.\newline

The following proposition shows that ``$\epsilon $-contraction'' is
equivalent to ''a big ratio between the first and second diagonal terms in
the $KAK$ decomposition.''

\begin{prop}
\label{g contracts}Let $\epsilon <\frac{1}{4}$. If $|\frac{a_{2}(g)}{a_{1}(g)%
}|\leq \epsilon ^{2}$, then $[g]$ is $\epsilon $-contracting. More
precisely, writing $g=k_{g}a_{g}k_{g}^{\prime }$, one can take $H_{g}$ to be
the projective hyperplane spanned by $\{k^{\prime
}{}_{g}^{-1}(e_{i})\}_{i=2}^{n}$, and $v_{g}=[k_{g}(e_{1})]$.

Conversely, suppose $[g]$ is $\epsilon $-contracting and $k$ is
non-archimedean with uniformizer $\pi $ (resp. archimedean), then $|\frac{%
a_{2}(g)}{a_{1}(g)}|\leq \frac{\epsilon ^{2}}{\left| \pi \right| }$ (resp. $|%
\frac{a_{2}(g)}{a_{1}(g)}|\leq 4\epsilon ^{2}$).
\end{prop}

\proof%
Suppose that $|\frac{a_{2}(g)}{a_{1}(g)}|\leq \epsilon ^{2}$. Using the
Cartan decomposition and the fact that $K$ acts by isometries, we can again
assume that 
\begin{equation*}
g=a_{g}=diag(a_{1},...,a_{n})
\end{equation*}
Let $[v]\in \Bbb{P}^{n-1}(k)$ be outside the $\epsilon $-neighborhood of $%
H_{g}$, which, in our case, simply means that (cf. (\ref{lin})) 
\begin{equation*}
d([v],H_{g})=\frac{\left| v_{1}\right| }{\left\| v\right\| }\geq \epsilon
\end{equation*}
Then 
\begin{equation*}
d([gv],[e_{1}])=\frac{\left\| gv\land e_{1}\right\| }{\left\| gv\right\| }%
\leq \frac{|a_2|\cdot\|v\|}{|a_1|\cdot|v_1|} \leq \epsilon
\end{equation*}
since $\left\| gv\land e_{1}\right\| \leq \left| a_{2}\right| \left\|
v\right\| $ and $\left\| gv\right\| \geq \left| a_{1}v_{1}\right| $.

The converse is more delicate. We shall describe the non-archimedean case,
and remark that the archimedean case can be dealt with in an analogous way. Again
we can take $g=a_{g}$. Let $f=(f_{1},...,f_{n})$ be a linear form of norm $1$
such that $\ker (f)=H_{g}$ and let $w\in k^n$ be a normalized representative
of the attracting point $v_g$. The fact that $[g]$ is $\epsilon $%
-contracting means that for every non-zero vector $v\in k^{n}$%
\begin{equation}
\left| f(v)\right| \geq \epsilon \left\| v\right\| \Rightarrow \left\|
gv\land w\right\| \leq \epsilon \left\| gv\right\| .  \label{cont}
\end{equation}
Suppose that for some index $i_{0}$, $\left| f_{i_{0}}\right| \geq \epsilon $%
. Take $v=e_{i_0}$ and get that 
\begin{equation*}
\left| a_{i_{0}}\right| \left\| e_{i_{0}}\land w\right\| \leq \epsilon
\left| a_{i_{0}}\right|
\end{equation*}
So $\left\| e_{i_{0}}\land w\right\| =\max_{i\neq i_{0}}\left| w_{i}\right|
\leq \epsilon $ and hence $\left| w_{i_{0}}\right| =1$, and $i_{0}$ is
unique. Since $\left\| f\right\| =\max \left| f_{i}\right| =1$, we must then
have $\left| f_{i_{0}}\right| =1$ and $\max_{i\neq i_{0}}\left| f_{i}\right|
<\epsilon $.

We now claim that $i_{0}=1$. Suppose the contrary and take $v=e_{1}+e_{i_{0}}
$. Then $\left\| v\right\| =1$ and $f(v)=f_{1}+f_{i_{0}}$. Since $\left|
f_{i_{0}}\right| =1\geq \epsilon >\left| f_{1}\right| $ we indeed have 
\begin{equation*}
\left| f(v)\right| \geq \epsilon \left\| v\right\| 
\end{equation*}
Hence, by (\ref{cont}), we obtain 
\begin{equation*}
\left\| a_{1}e_{1}\land w+a_{i_{0}}e_{i_{0}}\land w\right\| \leq \epsilon
\left\| a_{1}e_{1}+a_{i_{0}}e_{i_{0}}\right\| 
\end{equation*}
which also reads 
\begin{equation*}
\left| a_{1}\right| \leq \epsilon \left| a_{1}\right| 
\end{equation*}
and gives the desired contradiction.

Finally take $v=xe_{1}+e_{2}$, where $x\in k$ is chosen such that $\left|
x\right| $ is least possible and $\geq \epsilon $. Then $\left\| v\right\|
=1 $ and $f(v)=xf_{1}+f_{2}$. Since $\left| xf_{1}\right| \geq \epsilon
>\left| f_{2}\right| $ we have 
\begin{equation*}
\left| f(v)\right| \geq \epsilon \left\| v\right\|
\end{equation*}
and, again, by (\ref{cont}), we obtain 
\begin{equation*}
\left\| xa_{1}e_{1}\land w+a_{2}e_{2}\land w\right\| \leq \epsilon \left\|
xa_{1}e_{1}+a_{2}e_{2}\right\|
\end{equation*}
Suppose $\left| a_{2}\right| >\epsilon \left| x\right| \left| a_{1}\right| $
then the last inequality translates to 
\begin{equation*}
\left| a_{2}\right| \leq \epsilon \max \{\left| xa_{1}\right| ,\left|
a_{2}\right| \}
\end{equation*}
which leads to either $\left| a_{2}\right| \leq \epsilon \left| a_{2}\right| 
$, which is absurd, or $\left| a_{2}\right| \leq \epsilon \left|
x\right|\cdot \left| a_{1}\right| $, which contradicts the assumption.

Therefore we have obtained that $\left| a_{2}\right| \leq \epsilon \left|
x\right| \left| a_{1}\right| $, which, considering the choice of $x,$
implies the desired conclusion: 
\begin{equation*}
\left| \frac{a_{2}}{a_{1}}\right| \leq \frac{\epsilon ^{2}}{\left| \pi
\right|}
\end{equation*}
\endproof%

Note that the factor $\left| \pi \right| $ is not necessary when $\epsilon $
belongs to the value group of $k$.\newline

\begin{lem}
\label{open}Let $r,\epsilon \in (0,1]$. If $|\frac{a_{2}(g)}{a_{1}(g)}|\leq
\epsilon ,$ then $[g]$ is $\frac{\epsilon }{r^{2}}$-Lipschitz outside the $r$%
-neighborhood of the repulsive hyperplane $H_{g}=[span\{k^{\prime
}{}_{g}^{-1}(e_{i})\}_{i=2}^{n}]$.
\end{lem}

\proof%
Again, we can assume $g=a_{g}$, $v_{g}=[e_{1}]$, $H_{g}=[span(e_{i})_{i\geq
2}].$ Let $v,w$ be arbitrary non-zero vectors in $k^{n}$. Then $\left\|
gv\land gw\right\| \leq \left| a_{1}a_{2}\right| \left\| v\land w\right\| $
and $\left\| gv\right\| \geq \left| a_{1}v_{1}\right| $ and $\left\|
gw\right\| \geq \left| a_{1}w_{1}\right| $. Therefore 
\begin{equation}
d\big( [gv],[gw]\big)\leq \frac{|a_{2}|\cdot \left\| v\right\|\cdot \left\|
w\right\| }{|a_1|\cdot\left| v_{1}\right|\cdot \left| w_{1}\right| } d\big( %
[v],[w]\big).  \label{lipp}
\end{equation}
We conclude by observing that the $r$-neighborhood of $H_{g}=[span(e_{i})_{i%
\geq 2}]$ corresponds to non-zero vectors $v$ for which $\left| v_{1}\right|
\leq r\left\| v\right\| $. 
\endproof%

Note in particular that if $\delta >0,$ any projective transformation $%
[g]\in $PSL$_{n}(k)$ is $(1+\delta )$-Lipschitz in some open set of $\Bbb{P}%
^{n-1}(k)$.\newline

The following lemma gives a converse statement to the last result as well as
a handy criterion for $\epsilon $-contraction.

\begin{lem}
\label{a1/a2}Let $[g]$ be a projective transformation corresponding to $g\in 
\text{SL}_{n}(k)$, and assume that the restriction of $[g]$ to some open set 
$O\subset {\mathbb{P}}^{n-1}(k)$ is $\epsilon $-Lipschitz for some $\epsilon
<1$, then $|\frac{a_{2}(g)}{a_{1}(g)}|\leq \epsilon $ when $k$ is
non-archimedean, and $|\frac{a_{2}(g)}{a_{1}(g)}|\leq \frac{\epsilon }{\sqrt{%
1-\epsilon ^{2}}}$ when $k$ is archimedean.
\end{lem}

\proof%
Once again, we may assume $g=a_{g}$. Let $v$ be a non-zero vector in $k^{n}$
such that $[v]\in O$. For some $\delta \in k\backslash \{0\}$ small enough,
both $[w_{1}]=[v+\delta e_{1}]$ and $[w_{2}]=[v+\delta e_{2}]$ belong to $O$%
. Hence, for $w_{1}$%
\begin{equation}
\left| \delta a_{1}\right| \frac{\left\| av\land e_{1}\right\| }{\left\|
av\right\| \cdot\left\| aw_{1}\right\| }=d\big( [av],[aw_{1}]\big)\leq
\epsilon \cdot d\big( [v],[w_{1}]\big) =\epsilon \left| \delta \right| \frac{%
\left\| v\land e_{1}\right\| }{\left\| v\right\|\cdot \left\| w_{1}\right\| 
}  \label{ine}
\end{equation}
thus 
\begin{equation*}
\frac{\left\| av\land e_{1}\right\| }{\left\| av\right\| }\leq \epsilon
\end{equation*}
In the non-archimedean case, this implies $\left| a_{1}v_{1}\right| =\left\|
av\right\| $, and in the archimedean case, $\left\| av\land e_{2}\right\|
\geq \left| a_{1}v_{1}\right| \geq \sqrt{1-\epsilon ^{2}}\left\| av\right\| $%
. Now, expressing the Lipschitz condition for $v$ and $w_{2}$ as in (\ref{ine})
yields 
\begin{equation*}
\left| a_{2}\right| \frac{\left\| av\land e_{2}\right\| }{\left\|
av\right\| }\leq \epsilon \left| a_{1}\right|
\end{equation*}
which, in the non-archimedean case gives $\left| a_{2}/a_{1}\right| \leq
\epsilon $ and in the archimedean case, $\left| a_{2}/a_{1}\right| \leq
\epsilon /\sqrt{1-\epsilon ^{2}}$.%
\endproof%

\begin{defn}
A projective transformation $[g]\in $PSL$_{n}(k)$ is called \textbf{$%
(r,\epsilon )$-proximal} ($r>2\epsilon >0$) if it is $\epsilon $-contracting
with respect to some attracting point $v_{g}\in \Bbb{P}^{n-1}(k)$ and some
repulsive hyperplane $H_{g}$, such that $d(v_{g},H_{g})\geq r$. The
transformation $[g]$ is called \textbf{$(r,\epsilon )$-very proximal} if
both $[g]$ and $[g]^{-1}$ are $(r,\epsilon )$-proximal.
\end{defn}

Similar notions of $(r,\epsilon )$-proximality were defined in \cite{abels}
and \cite{ben}.

\begin{defn}
\label{sep}A \textit{finite} subset $F\subset \text{PSL}_{n}(k)$ is called $%
(m,r\mathbf{)}$\textbf{-separating} ($r>0$, $m\in \Bbb{N}$) if for every
choice of $2m$ points $v_{1},...,v_{2m}$ in $\Bbb{P}^{n-1}(k)$ and $2m$
projective hyperplanes $H_{1},...,H_{2m}$ there exists $\gamma \in F$ such
that 
\begin{equation*}
\min_{1\leq i,j\leq 2m}\{d(\gamma v_{i},H_{j}),d(\gamma
^{-1}v_{i},H_{j})\}>r.
\end{equation*}
\end{defn}

We use the properties described above in order to construct $\epsilon $-very
contracting and $(r,\epsilon )$-very proximal elements from two basic
ingredients : an $\epsilon $-contracting element and an $r$-separating set.

\begin{prop}
\label{P-VP}Let $F$ be a $(1,r)$-separating set ($r<1$) in $\text{PSL}_{n}(k)
$ with uniform bi-Lipschitz constant $C=\max \{\text{biLip}\big( [f]\big) %
:[f]\in F\}$ (see lemma \ref{biLip}).

$(i)$ If $\epsilon <\frac{r}{2C}$ and $[g]\in $PSL$_{n}(k)$ is an $\epsilon $%
-contracting transformation, then one can find $[f]\in F$ such that $[fg]$
is $(r,C\epsilon )$-proximal.

$(ii)$ Let $d=4$ when $k$ is archimedean and $d=\frac{1}{\left| \pi \right| }
$ when $k$ is non-archimedean. If $\epsilon <\frac{r}{\sqrt{2Cd}}$ and $%
[g]\in $PSL$_{n}(k)$ is $\epsilon $-contracting, then there is an element $%
[f]\in F$, such that

\begin{equation*}
\lbrack gfg^{-1}]
\end{equation*}
is $\frac{\sqrt{2Cd}}{r}\epsilon $-very contracting.

$(iii)$ If $\epsilon <\frac{r}{2C^{2}}$ and $[g]\in $PSL$_{n}(k)$ is an $%
\epsilon $-very contracting transformation, then there is an element $[f]\in F$,
such that $[gf]$ is $(\frac{r}{C},C\epsilon )$-very proximal.
\end{prop}

\proof%
Let $v_{g}$ (resp. $H_{g}$) be an attracting point (resp. repulsive
hyperplane) corresponding to the transformation $[g]$. Let $[f]\in F$ be
such that $d([f]v_{g},H_{g})>r$. Then $[fg]$ is clearly $C\epsilon $%
-contracting with attracting point $[f]v_{g}$ and repulsive hyperplane $%
H_{g}.$ The choice of $f$ guarantees that $[fg]$ is $(r,C\epsilon )$%
-proximal. This proves (i).

By lemma \ref{open} and the remark following it, there is an open set $O$ on
which $[g^{-1}]$ is, say, $\sqrt{2}$-Lipschitz. Let $[u]$ be some point
inside $O$. We can find $[f]\in F$ such that $d([fg^{-1}u],H_{g})>r$ and $%
d([f^{-1}g^{-1}u],H_{g})>r$. Since $[g]$ is $\epsilon $-contracting, by
proposition \ref{g contracts} we have $\left| \frac{a_{2}(g)}{a_{1}(g)}%
\right| \leq d\epsilon ^{2}$. Hence by lemma \ref{open}, $[g]$ is $\frac{%
d\epsilon ^{2}}{r^{2}}$-Lipschitz outside the $r$-neighborhood of $H_{g}$.
Therefore both $[gfg^{-1}]$ and $[gf^{-1}g^{-1}]$ are $\frac{dC\sqrt{2}%
\epsilon ^{2}}{r^{2}}$--Lipschitz in some open neighborhood of $[u]$, which
implies by lemma \ref{a1/a2}, that they both satisfy $\frac{|a_{2}|}{|a_{1}|}
\leq 2\frac{dC\epsilon ^{2}}{r^{2}}$. By proposition \ref{g contracts} it
now follows that they are both $\frac{\sqrt{2Cd}}{r}\epsilon $-contracting
transformations. This proves (ii).

Find $[f]\in F$ which takes $v_{g}$ at least $r$ away from $H_{g}$ and for which 
$[f^{-1}]$ takes $v_{g^{-1}}$ at least $r$ away from $H_{g^{-1}}$. Then $%
[gf] $ is $(\frac{r}{C},C\epsilon )$-proximal and $[f^{-1}g^{-1}]$ is $%
(r,C\epsilon )$-proximal. Hence $[gf]$ is $(\frac{r}{C},C\epsilon )$-very
proximal. This proves (iii).

\endproof%

We end this section by showing how to obtain generators of a free group via
the so-called ping-pong lemma, once we are given a very contracting element
and a separating set. The following definition and lemma are classical (cf. 
\cite{T}).

\begin{defn}[Ping-pong pair]
Let $k$ be a local field and $V$ a finite dimensional $k$-vector space. A
pair of projective transformations $a,b\in $PSL$(V)$ is called a \textbf{%
ping-pong pair} if both $a$ and $b$ are $(r,\epsilon )$-very proximal, with
respect to some $r>2\epsilon >0$, and if the attracting points of $a$ and $%
a^{-1}$ (resp. of $b$ and $b^{-1}$) are at least $r$-apart from the
repulsive hyperplanes of $b$ and $b^{-1}$ (resp. of $a$ and $a^{-1}$).

More generally, a $m$-tuple of projective transformations $a_{1},\ldots
,a_{m}$ is called a \textbf{ping-pong $m$-tuple} if all $a_{i}$'s are $%
(r,\epsilon )$-very proximal (for some $r>2\epsilon >0$) and the attracting
points of $a_{i}$ and $a_{i}^{-1}$ are at least $r$-apart from the repulsive
hyperplanes of $a_{j}$ and $a_{j}^{-1}$, for any $i\neq j$.
\end{defn}

We have the following variant of the ping-pong lemma (see \cite{T} 1.1).

\begin{lem}
If $a_{1},\ldots ,a_{m}\in $PSL$(V)$ form a ping-pong $m$-tuple, then they are free generators of a free group $F_m$.
\end{lem}

The following proposition explains how the above ingredients can be combined to yield a free group.

\begin{prop}
\label{free}Given $a_{1},\ldots ,a_{m}\in \text{PSL}(V)$, an $(m,r)$%
-separating set $F$, and an $\epsilon $-very contracting element $\gamma $
(with $\epsilon <\frac{r}{2c^{4}}$, where $c$ is the maximal bi-Lipschitz
constant of $\{a_{1},...,a_{m}\}\cup F$), there are $h_{1},\ldots ,h_{m}\in F
$ and $g_{2},\ldots ,g_{m}\in F$ such that 
\begin{equation*}
\gamma a_{1}h_{1},~g_{2}\gamma a_{2}h_{2},\ldots ,~g_{m}\gamma a_{m}h_{m}
\end{equation*}
form a ping-pong $m$-tuple of $(\frac{r}{c},c^{3}\epsilon )$-very proximal
transformations and hence generate a free group $F_m$.
\end{prop}

\proof
In this proof, whenever we speak about an $\epsilon $-contracting or $%
(r,\epsilon )$-proximal transformation, we shall choose and fix a point and
a projective hyperplane with respect to which it is $\epsilon $-contracting
or $(r,\epsilon )$-proximal and shall refer to them as the attracting point
and the repulsive hyperplane of that transformation.

As $\gamma $ is $\epsilon $-very contracting, $\gamma a_{1}$ is $c\epsilon $%
-very contracting. Now by proposition \ref{P-VP} $(iii)$, we can find $%
h_{1}\in F$ such that $\gamma a_{1}h_{1}$ is $(\frac{r}{c},c^{2}\epsilon )$%
-very proximal.

We proceed by induction and suppose that $x_{1}=\gamma_{1}a_{1}h_{1}$ and
$x_{j}=g_{j}\gamma a_{j}h_{j}$ have been constructed for all $%
j<i$ ($i=2,...,m$). Let's construct $h_{i}$ and $g_{i}$. First note that $%
\gamma a_{i}$ is $c\epsilon $-very contracting. We take $h_{i}\in F$ such
that $h_{i}^{-1}$ maps the attracting point of $(\gamma a_{i})^{-1}$ at a
distance at least $r$ from the repulsive hyperplanes of all transformations $%
x_{j}$ and $x_{j}^{-1}$ for $j<i$, and $h_i$ maps the attracting point of any
$x_j$ and $x_j^{-1}$ $(j<i)$ at least $r$ away from the repulsive hyperplane
of $\gc a_i$. 
This is possible since $F$ is an $(m,r)$-separating set. 
The second requirement for $h_i$ implies that the repulsive hyperplane
$H_{(\gc a_ih_i)}=h_i^{-1}(H_{(\gc a_i)})$ of $\gc a_ih_i$ is at least $r/c$ 
away from the attracting points of the $x_j$'s and $x_j^{-1}$'s.

Then pick $g_{i}\in F$ which takes the attracting point of $%
\gamma a_{i}h_i$ at a distance at least $r$ from the repulsive hyperplanes of
all $x_{j}$'s and $x_{j}^{-1}$'s for all $j<i$ and also from the repulsive
hyperplane of $\gamma a_{i}h_{i}$, and whose inverse $g_i^{-1}$ takes the attracting 
point of any $x_j$ and $x_j^{-1}$ $(j<i)$ at least $r$ away from
the repulsive hyperplane of $(\gamma a_{i}h_i)^{-1}$. This means that the 
repulsive hyperplane $g_i(H_{(\gc a_ih_i)^{-1}})=H_{(g_i\gc a_ih_i)^{-1}}$
of $g_i\gc a_ih_i$
is at least $r/c$ away from the attracting points all $x_j$ and $x_j^{-1}$.
Additionally we can require that $%
g_{i}^{-1}$ takes the attracting point of $(\gamma a_{i}h_{i})^{-1}$ at a
distance at least $r$ from the repulsive hyperplane of 
$(\gamma a_{i}h_i)^{-1}$. 

Then clearly, $x_{i}=g_{i}\gamma a_{i}h_{i}$ is $(\frac{r}{c}%
,c^{3}\epsilon )$-very proximal and $(x_{1},...,x_{i})$ form a ping-pong $i$%
-tuple of $(\frac{r}{c},c^{3}\epsilon )$-very proximal transformations. 
\endproof%


\section{Constructing very proximal elements in $\Gamma \label{cons}$}

Our aim in this section is to establish an action of $\Gamma $, on some
projective space over a local field, which has many very proximal elements.
Proposition \ref{P-VP} tells us how to get very proximal elements out of
contracting elements when a separating set is available. We shall show how
to construct an action of $\Gamma$, on a projective space over a local
field, with both ingredients : contracting elements and separating set.

For an algebraic number field there are naturally associated local fields -
its completions. The following lemma reduces the general case to this.

\begin{lem}
\label{alg}Let $\Bbb{G}$ be a semisimple algebraic group defined over $\Bbb{Q%
}$, and $G=\Bbb{G}(\Bbb{R})^{0}$ be the corresponding connected semisimple
Lie group. Let $\Gamma $ be a finitely generated dense subgroup in $G$.
Then, there exists a number field $K\subset {\mathbb{R}}$ and a group
homomorphism $\pi :\Gamma \to \Bbb{G}$ such that $\pi (\Gamma
)\subset \Bbb{G(}K)$ and $\pi (\Gamma )$ is still dense in $G$.
\end{lem}

\proof
Choose a finite set of generators $\gamma _{1},\gamma _{2},\ldots ,\gamma
_{k}$ such that $\gamma _{1},\gamma _{2}$ are close to $1$ and $(\gamma
_{1},\gamma _{2})$ belongs to the complement of the proper exponential
algebraic set introduced in theorem \ref{denseSS}. Let $\Gamma =\{\gamma
_{1},...,\gamma _{k}$ $|$ $(r_{\alpha })_{\alpha }\}$ be a presentation of $%
\Gamma $. The variety of representations 
\begin{equation*}
\text{Hom}\big(\Gamma ,\Bbb{G}\big) =\{(g_{1},...,g_{k})\in \Bbb{G}%
^{k}:\forall \alpha ,~r_{\alpha }(g_{1},...,g_{k})=1\}
\end{equation*}
is an algebraic subvariety of $\Bbb{G}^{k}$ which is defined over $\Bbb{Q}$
. Its set of algebraic points is therefore dense (for the real topology) in
its set of real points. Choosing an algebraic point in $\text{Hom}\big(%
\Gamma ,\Bbb{G(R)}\big)$ very close to the original $(\gamma _{1},...,\gamma
_{k})$, we obtain a representation $\pi $ of $\Gamma $ into $G$ such that $%
\pi (\Gamma )$ is contained in $\Bbb{G(}\overline{\Bbb{Q}})$ and is still
dense in $G$, as small deformations of $\gamma _{1},\gamma _{2}$ generate
dense subgroups. Finally, since it is finitely generated, it is contained in 
$\Bbb{G}(K)$ for some number field $K$.
\endproof

As explained in the last section, ``$\varepsilon $-contraction'' is
equivalent to ``big ratio between the first and the second diagonal entries
in the $KAK$ decomposition.''

\begin{lem}[Constructing contracting elements]
\label{proximal}Let $K$ be a number field and $\Bbb{G}$ be a non-trivial
semisimple Zariski-connected algebraic group defined over $K.$ Let $R$ be a
finitely generated subring of $K$ and $I$ an infinite subset of $\Bbb{G}(R)$%
. Then for some completion $k$ of $K$ there is an irreducible rational
representation $\rho :\Bbb{G\to SL}(V)$ defined over $k$, with $\dim
_{k}V\geq 2$, such that the set 
\begin{equation*}
\left\{ \frac{|a_{1}\big(\rho (g)\big)|}{|a_{2}\big(\rho (g)\big)|}:g\in
I\right\} 
\end{equation*}
is unbounded, where $a_{1}\big(\rho (g)\big)$ and $a_{2}\big(\rho (g)\big)$
are the first and second diagonal terms in a Cartan decomposition of $\rho
(g)$ in $\Bbb{SL}(V)$.
\end{lem}

\proof
As $R$ is finitely generated, the discrete diagonal embedding of $K$ in its
ad\`{e}le group gives rise to a discrete embedding of $R$ into a product of
finitely many places. Explicitly, there is a finite set $S$ of places of $K,$
including all archimedean ones, such that $R$ is contained in the ring of $S$%
-integers $\mathcal{O}_{K}(S)$. Projecting to the finite product of places
corresponding to $S$, we get a discrete embedding 
\begin{equation*}
R\hookrightarrow \prod_{\nu \in S}K_{\nu },
\end{equation*}
which gives rise to a discrete embedding 
\begin{equation*}
\Bbb{G}(R)\hookrightarrow \prod_{\nu \in S}\Bbb{G}(K_{\nu }).
\end{equation*}
The infinite set $I$ is mapped to some unbounded set in the latter product.
Hence for some place $\nu \in S$, the embedding $\Bbb{G}(R)\hookrightarrow 
\Bbb{G}(K_{\nu })$ sends $I$ to an unbounded set. Take $k=K_{v}$. Since $%
\Bbb{G}$ is also defined over $k$, there exists a $k$-rational faithful
representation of $\Bbb{G}$ on some $k$-vector space $V_{0}.$ Under this
representation, the set $I$ is sent to some unbounded set in $\Bbb{SL(}V_{0})
$ which we shall continue to denote by $I$.

Fix a $k$-basis of $V_{0}$ and consider the Cartan decomposition in $\Bbb{SL(%
}V_{0})$ corresponding to this basis. Recall the notations of section \ref
{dyn}. For $x\in \Bbb{SL(}V_{0}),~|a_{1}(x)|\geq |a_{2}(x)|\geq \ldots \geq
|a_{d}(x)|>0$ denote the corresponding diagonal coefficients in the Cartan
decomposition of $x$ in $\Bbb{SL(}V_{0})$. Since their product is $1$, there
is some $i~(1\leq i<d)$ for which the set $\{\frac{|a_{i}(g)|}{|a_{i+1}(g)|}\}_{g\in
I}$ is unbounded.

Considering the $i$'th wedge product representation of the above
representation of $\Bbb{G,}$ we obtain a rational representation 
\begin{equation*}
\rho _{0}:\Bbb{G\to SL(}\bigwedge^{i}V_{0})
\end{equation*}
defined over $k$. Looking at the Cartan decomposition in $\Bbb{SL(}%
\bigwedge^{i}V_{0})$ with respect to the basis of $\bigwedge^{i}V_{0}$
induced by the given basis of $V_{0},$ we see that $a_{1}\big(\rho _{0}(x)%
\big) =a_{1}(x)...a_{i}(x)$ and $a_{2}\big(\rho _{0}(g)\big) %
=a_{1}(x)...a_{i-1}(x)a_{i+1}(x)$, and hence, under this representation, the
set 
\begin{equation*}
\left\{\frac{|a_{1}\big(\rho _{0}(g)\big)|}{|a_{2}\big(\rho _{0}(g)\big)|}%
\right\}_{g\in I}
\end{equation*}
is unbounded.

The representation $\rho _{0}$ might be reducible. However, as $\Bbb{G}(k)$
is a connected semisimple algebraic group, the $k$-rational representation $%
\rho _{0}$ is completely reducible. Let $W=\bigwedge\nolimits^{i}V_{0}$ be
the representation space. Decompose it into $\Bbb{G}(k)$-irreducible spaces $%
W=W_{1}\oplus ...\oplus W_{q}$. Any $x=\rho _{0}(g)\in \rho _{0}\big(\Bbb{G}%
(k)\big)\leq \Bbb{SL(}W\Bbb{)}$ stabilizes each $W_{i}$, and we shall write $%
x_{i}\in \Bbb{SL}(W_{i})$ for its restriction to the subspace $W_{i}$ (note
that the determinant of $x_{i}$ is 1 as $\Bbb{G}(k)$ is semisimple).

The space $W$ is endowed with the norm $\left\| \cdot \right\| $
corresponding to the above choice of a basis of $V_{0}$ (as in section \ref
{dyn}) for which $|a_{1}(g)|=\left\| g\right\| $ and $|a_{1}(g)a_{2}(g)|=\left\|
\bigwedge^{2}g\right\| $ for any $g\in \Bbb{SL}(W)$. Similarly we can choose
a basis and corresponding norms $\left\| \cdot \right\| _{i}$ for each $%
W_{i}$ such that 
\begin{equation*}
|a_{1}(x_{i})|=\left\| x_{i}\right\|_{i} \text{~and~}
|a_{1}(x_{i})a_{2}(x_{i})|=\left\| \bigwedge^{2}x_{i}\right\| _{i}.
\end{equation*}
As any two norms on the finite dimensional vector space $M_n(k),~\big( %
n=\dim (W)\big)$ are equivalent, there is $c_{1}>0$ with respect to which 
\begin{equation*}
\frac{1}{c_{1}}\max_{1\leq j\leq q}\{\left\| x_{j}\right\| _{j}\}\leq
\left\| x\right\| \leq c_{1}\max_{1\leq j\leq q}\{\left\| x_{j}\right\|
_{j}\}.
\end{equation*}
Similarly, for some constant $c_2>0$ 
\begin{equation*}
\frac{1}{c_{2}}\max_{i,j}\{\left\|
(\bigwedge\nolimits^{2}x)_{ij}\right\|_{ij} \}\leq \left\|
\bigwedge\nolimits^{2}x\right\| \leq c_{2}\max_{i,j}\{\left\|
(\bigwedge\nolimits^{2}x)_{ij}\right\| _{ij}\}.
\end{equation*}

Replacing $I$ by an infinite subset if necessary, we can find an index $%
i_{0}\leq q$ with $\left\| x\right\| \leq c_{1}\left\| x_{i_{0}}\right\|
_{i_{0}}$ for any $x$ of the form $x=\rho _{0}(g)$, $g\in I$. As $%
\{\left\|\rho _{0}(g)\right\|\}_{g\in I}$ is unbounded, and each $x_{i}$
belongs to $\Bbb{SL}(W_{i})$, we must have $\dim W_{i_{0}}\geq 2$. Moreover, we have 
\begin{equation*}
\left\|\bigwedge\nolimits^{2}x\right\| \geq \frac{1}{c_{2}}\left\|
(\bigwedge\nolimits^{2}x)_{ii}\right\| _{ii}=\frac{1}{c_{2}}\left\|
\bigwedge\nolimits^{2}x_{i}\right\| _{i},
\end{equation*}
which implies 
\begin{equation*}
\frac{|a_{1}(x)|}{|a_{2}(x)|}\leq c_{1}^{2}c_{2}\frac{|a_{1}(x_{i_{0}})|}{
|a_{2}(x_{i_{0}})|}.
\end{equation*}
Hence $\rho=(\rho _{0})_{|Wi_{0}}$ is the desired rational irreducible
representation.

\endproof

\begin{lem}[Constructing separating sets]
\label{sepa} Let $m$ be a positive integer, $k$ a local field, $\Bbb{G}$ a
Zariski-connected $k$-algebraic group, $V$ a $k$-vector space, and $\rho :%
\Bbb{G\to SL(}V\Bbb{)}$ a $k$-irreducible rational representation.
Then for any Zariski-dense subset $\Omega \subset \Bbb{G(}k)$, there exists
a finite subset $F\subset \Omega $ and a positive real number $r>0$, which
gives rise under $\rho $ to an $(m,r)$-separating set (cf. \ref{sep}) on $%
\Bbb{P}(V)$.
\end{lem}

\proof%
For each $\gamma \in \Gamma $, let $M_{\gamma }$ be the set of all tuples $%
(v_{i},H_{i})_{i=1}^{2m}$ of $2m$ points $v_{i}\in \Bbb{P}(V)$ and $2m$
projective hyperplanes $H_{i}\subset \Bbb{P}(V)$, not necessarily different,
such that for some $1\leq i,j\leq 2m$, $\gamma \cdot v_{i}\in H_{j}$ or $%
\gamma ^{-1}\cdot v_{i}\in H_{i}$.

We first claim that 
\begin{equation}
\bigcap_{\gamma \in \Omega }M_{\gamma }=\emptyset .  \label{0}
\end{equation}
Suppose this were not the case. Then there would exist $%
(v_{i},H_{i})_{i=1}^{2m}$ such that 
\begin{equation*}
\Omega \subset \bigcup_{1\leq i,j\leq 2m}\left\{ \gamma \in \Bbb{G}%
(k),\gamma \cdot v_{i}\in H_{j}\text{~or~}\gamma ^{-1}\cdot v_{i}\in
H_{j}\right\} .
\end{equation*}
The sets $\left\{ \gamma \in \Bbb{G}(k):\gamma \cdot v_{i}\in H_{j}\right\} $
and $\left\{ \gamma ^{-1}\in \Bbb{G}(k):\gamma \cdot v_{i}\in H_{j}\right\} $
are clearly Zariski-closed and proper, since $\Bbb{G}(k)$ acts irreducibly
on $V$. As $\Bbb{G}(k)$ is Zariski-connected, this yields a contradiction.
Thus we have $(\ref{0})$.

Second, as the sets $M_{\gamma }$ are compact in the appropriate product of
grassmannians, there is a finite subset $F\subset \Omega $, such that 
\begin{equation}
\bigcap_{\gamma \in F}M_{\gamma }=\emptyset .  \label{1}
\end{equation}

The following MaxMin of continuous functions 
\begin{equation}
\max_{\gamma \in F}\Big(\min_{1\leq i,j\leq 2m}\{d(\gamma \cdot
v_{i},H_{j}),d(\gamma ^{-1}\cdot v_{i},H_{j})\}\Big)  \label{2}
\end{equation}
depends continuously on $(v_{i},H_{i})_{i=1}^{2m}$, and by $(\ref{1})$ never
vanishes. Finally, the compactness of the set of all tuples $%
(v_{i},H_{i})_{i=1}^{2m}$ in 
\begin{equation*}
\big(\Bbb{P}(V)\times \Bbb{G}r_{\dim (V)-1}(V)\big)^{2m}
\end{equation*}
implies that $(\ref{2})$ attains a positive minimum $r$. Thus $F$ gives rise
to the desired $(m,r)$-separating set.%
\endproof%

Putting together proposition \ref{P-VP} with lemmas \ref{proximal} and \ref
{sepa}, we obtain a great liberty in the choice of very contracting (or very
proximal) elements in $\Gamma $.

\begin{lem}\label{mrr}
Let $m$ be a positive integer, $K$ a number field, $\Bbb{G}$ a
Zariski-connected semisimple $K$-algebraic group, $R\subset K$ a finitely
generated subring, $I\subset \Bbb{G}(R)$ an infinite subset, and $\Omega
\subset \Bbb{G(}K)$ a Zariski-dense subset of $\Bbb{G}$. Then there exist:

\begin{itemize}
\item  An embedding $\sigma :K\hookrightarrow k$ of $K$ into some local
field $k$.

\item  An irreducible rational representation $\rho :\Bbb{G}\to \Bbb{%
SL}(V)$ defined over $k$,

\item  A finite $(m,r)$-separating set $F\subset \Omega $ (for some $r>0$)
for the action induced by $\rho $ on the projective space $\Bbb{P}(V)$,
\end{itemize}

such that for any $\varepsilon >0$, there is $g\in I$ and $f\in F$ for which 
$gfg^{-1}$ acts as an $\varepsilon $-very contracting transformation on $%
\Bbb{P}(V)$.\newline
\end{lem}

\proof

Lemma $(\ref{proximal})$ yields the local field $k$, the embedding $%
K\hookrightarrow k$, and the representation $\rho $, while lemma $(\ref{sepa}%
)$ yields the $r$-separating set $F\subset\Omega$.

For any $\varepsilon >0$, we can find $g\in I$ with 
\begin{equation*}
\frac{a_{1}\big(\rho (g)\big)}{a_{2}\big(\rho (g)\big)}\geq \frac{r}{%
\epsilon \sqrt{2Cd}}
\end{equation*}
where $C=\max \{\text{biLip}(f):f\in F\}$ and $d$ is as in proposition $(\ref
{P-VP})$ ($=4$ in archimedean case, $\frac{1}{\left| \pi \right|}$ in the
non-archimedean case). By lemma $(\ref{g contracts})$, $g$ acts as a $\frac{%
\epsilon \sqrt{2Cd}}{r}$-contracting transformation on $\Bbb{P(}V)$.
Finally, for $\epsilon $ small enough, proposition $(\ref{P-VP})$ provides $%
f\in F$ such that $gfg^{-1}$ is $\epsilon $-very contracting.

\endproof

Finally, we obtain the following result:

\begin{thm}\label{main} 
Let $\Bbb{G}$ be a Zariski connected semisimple algebraic $K$-group, where $K$ is a number field. Let $R$ be a finitely generated subring of $K$, and $\gO\subset\Bbb{G}(R)$ a 
Zariski dense subset of $\Bbb{G}$ with $\gO =\gO^{-1}$.
Then for any $m\in\BN$ and any $a_1,\ldots ,a_m\in \Bbb{G}(K)$ there are $x_1,\ldots ,x_m$ with $x_i\in \gO^4a_i\gO$ which are generators of a free group $F_m$.
\end{thm}

\begin{proof}
This follows from the combination of lemma \ref{mrr} and proposition \ref{free}.
\end{proof}

\begin{rem}
Note that the original statement of Tits' alternative (\cite{T}, theorem 1)
can be easily derived from the last theorem, first by taking a
homomorphism of $\Gamma $ into $\Bbb{GL}_{n}(K)$, for some number field $K$,
whose image is not almost solvable (this is possible as the index of the
solvable subgroup of an almost solvable group and the length of the derived
series of solvable groups are both uniformly bounded for groups contained in 
$\Bbb{GL}_{n}$) and then projecting to a semisimple factor of the Zariski
closure. Then take $\Omega =\Gamma $ in the last theorem.
\end{rem}


\section{Proofs of theorems \ref{th1}, \ref{th2}, and \ref{th3}}

To complete the proofs of the main results of this paper, we shall need some
preliminary lemmas.

\begin{lemma}
\label{fin}\label{gen} Let $G$ be a connected real Lie group and $\Gamma
\subset G$ a dense subgroup. Then $\Gamma $ is generated by $\Gamma \cap U$,
for any identity neighborhood $U$ in $G$.
\end{lemma}

\proof%
Let $H$ be the subgroup of $G$ generated by $\Gamma \cap U$. Since $\Gamma $
is dense and $G$ connected, $H$ is again dense in $G$. If $\gamma \in \Gamma 
$ we can thus find $h\in H$ such that $h\in \gamma U$. But then $\gamma
^{-1}h\in \Gamma \cap U\subset H$. Hence $\Gamma =H$.%
\endproof%

\begin{lemma}
\label{Zar} Let $\Bbb{G}$ be a Zariski-connected algebraic group defined
over $\Bbb{R}$ and let $G=\Bbb{G(R)}^{0}$ be the corresponding connected
real Lie group. Let $\Gamma $ be a dense subgroup in $G$ and $\pi :\Gamma
\to G$ a group homomorphism such that $\pi (\Gamma )$ is
Zariski-dense. Then for any identity neighborhood $U\subset G$, $\pi (U\cap
\Gamma )$ is Zariski-dense in $\Bbb{G}$.
\end{lemma}

\proof%
Let $(U_{n})_{n}$, $U_{n}\supset U_{n+1}$, be a nested sequence of identity
neighborhoods in $G$, such that for any identity neighborhood $U$, $%
U_{n}\subset U$ for large $n$'s. Let $U_{n}^{\prime }=U_{n}\cap \Gamma $ and
let $X_{n}$ be the Zariski closure of $\pi (U_{n}^{\prime })$ in $G$. Since $%
(X_{n})_{n}$ is a decreasing sequence of Zariski-closed subsets of $G$, it
stabilizes at some finite step $n_{0}$. For any $g\in U_{n_{0}}^{\prime }$, $%
gU_{n}^{\prime }\subset U_{n_{0}}^{\prime }$, for $n$ large enough, hence $%
\pi (g)X_{n_{0}}=\pi (g)X_{n}\subset X_{n_{0}}$. By lemma \ref{gen}, $%
U_{n_{0}}^{\prime }$ generates $\Gamma $, hence $\pi (\Gamma
)X_{n_{0}}\subset X_{n_{0}}$. Since $\pi (\Gamma )$ is Zariski dense, this
implies $X_{n_{0}}=G$. Thus $\pi (U\cap \Gamma )$ is Zariski-dense for any
identity neighborhood $U\subset G$.%
\endproof%

\proof[\bf Proof of theorem \ref{th1}] The adjoint group $H=\text{Ad}(G)$ is
a center free connected semisimple Lie group, hence is of the form $\Bbb{H(R}%
)^{0}$ where $\Bbb{H}$ is some Zariski-connected semisimple $\Bbb{Q}$%
-algebraic group.

By Corollary \ref{01}, we can assume that $\Gamma $ is finitely generated.
By lemma \ref{alg}, there is a homomorphism $\pi $ of $\Gamma $ into $\text{%
Ad}(G)$ with dense, hence Zariski-dense, image which lies in $\Bbb{H}(K)$
for some number field $K$. Since $\Gamma $ is finitely generated, there is a
finitely generated subring $R$ of $K$, such that $\pi (\Gamma )\subset \Bbb{H%
}(R)$.

Consider two elements $a$ and $b$ in $\Gamma $ and a sufficiently small identity
neighborhood $U^{\prime }$ in $G$ such that for any $x\in V(a):=U^{\prime
}aU^{\prime }$ and $y\in V(b):=U^{\prime }bU^{\prime }$, $\langle x,y\rangle 
$ is a dense subgroup of $G$. By theorem \ref{denseSS} this is always
possible. Let $U$ be a smaller symmetric identity neighborhood in $G$ such that $%
U^{4}\subset U^{\prime }$ and denote 
\begin{equation*}
\Omega =\pi (\Gamma \cap U).
\end{equation*}
By lemma \ref{Zar}, the conditions of theorem \ref{main} are fulfilled.
We thus obtain the desired dense free subgroup.

\endproof%

The proof of theorem \ref{th2} is basically the same as the proof of theorem 
\ref{th1} and we shall not go over it again, but only describe the needed
modifications in the above argument. \newline
\newline
\proof[\bf Proof of theorem \ref{th2} (outline)]

Pick $a_{1},\ldots ,a_{m}\in \Gamma $ close to the identity, and a small
symmetric identity neighborhood $U$ such that for any selection $x_{i}\in
U^{4}a_{i}U~,i=1,\ldots ,m$, the group $\langle x_{1},\ldots
,x_{m}\rangle $ is dense in $G$. The existence of such $a_{i}$'s and $U$ is
an easy consequence of theorem \ref{denseTP}.

Now project $G$ on the semisimple quotient $S=G/R$ where $R$ is the radical
of $G,$ and find as above a projective space over a local field and an action of $%
\Gamma $ (via $S$) for which there is

\begin{itemize}
\item  a $(m,r)$-separating set $F\subset U\cap \Gamma $ (for some $r>0$),
and

\item  for $\epsilon <\frac{r}{2c^{4}},$ (where $c$ is the maximal
bi-Lipschitz constant of $\{a_{1},...,a_{m}\}\cup F$ when acting on the
projective space), an element  $\gamma _{\epsilon }\subset U^{3}\cap \Gamma $ acting as
an $\epsilon $-very contracting element.
\end{itemize}

Then there are 
\begin{equation*}
h_{i}\in F~(1\leq i\leq m)\text{ and }g_{i}\in F~(2\leq i\leq m),
\end{equation*}
such that 
\begin{equation*}
\langle \gamma _{\epsilon }a_{1}h_{1},~g_{2}\gamma _{\epsilon
}a_{2}h_{2},\ldots ,~g_{m}\gamma _{\epsilon }a_{m}h_{m}\rangle
\end{equation*}
is dense and free. \qed

\begin{proof}[\bf Proof of theorem \ref{th3}]
Recall that $G$ is a connected non-solvable Lie group and $\gC\leq G$ is a 
finitely generated dense subgroup. 
Define $G_0=G$ and $G_n=\overline{[G_{n-1},G_{n-1}]}$. Then as 
$\dim (G) <\infty$ the sequence $G_n$ stabilizes at some finite step $k$,
$H=G_k$ is a connected topologically perfect closed normal subgroup of 
$G$, and the quotient $G/H$ is solvable. Moreover it follows by induction
on $n$ that $\gC\cap G_n$ is dense in $G_n$.

Let $l$ be the codimension of $H$ and let $m$ be the minimal number of 
generators of the Lie algebra of $H$, then $l+m\leq\dim (G)$.
Since the projection of $\gC$ to $G/H$ is finitely generated and dense,
proposition \ref{denseInFG} implies that it contains $2l$ elements which 
generate a dense subgroup in $G/H$. Let $a_1,\ldots ,a_{2l}\in\gC$ be 
arbitrary lifts of these elements.
Additionally, as in the proof of theorem \ref{th2},
let $a_{2l+1},\ldots ,a_{2l+m}$ be $m$ elements in
$H\cap\gC$ near the identity, and let $U\subset H$ be a small identity 
neighborhood of $H$, such that for any selection 
$x_i\in U^4a_iU~(2l+1\leq i\leq 2l+m)$, the 
group $\langle x_{2l+1},\ldots ,x_{2l+m}\rangle$ is dense in $H$. Then 
for any selection $x_i\in U^4a_iU~(1\leq i\leq 2l+m)$ the group
$\langle x_1,\ldots ,x_{2l+m}\rangle$ is dense in $G$.

From this point, we can continue exactly as in the proof of theorem \ref{th2}
in order to find $x_i$'s such that $\langle x_1,\ldots ,x_{2l+m}\rangle$
is free.
\end{proof}

\begin{remark}
The bound $2\dim (G)$ is not sharp as the latter proof shows. But $2l+m$ can
also be decreased in some cases. For example, it is easy to see that if $G$
is a direct product of a semisimple group with ${\mathbb{R}}^{d}$, then $%
\max \{2,2d\}$ is the sharp bound for theorem \ref{th3}.\newline
\end{remark}

\begin{Ack}
We would like to thank A. Lubotzky, from whom we learned about the question,
as well as G.A. Margulis for their help and support while working on this
paper. We are also thankful to Y. Benoist and S. Mozes for many valuable
discussions and comments and to Y. Guivarc'h for pointing out to us the work
of Carri\`{e}re and Ghys and the link with the theory of amenable actions.
\end{Ack}



\begin{thebibliography}{99}
\bibitem{abels}  Abels, H.; Margulis, G. A.; Soifer, G. A. Semigroups
containing proximal linear maps, \textit{Israel J. Math.} \textbf{91}
(1995), no. 1-3, 1--30.

\bibitem{ben}  Benoist, Y. Propri\'{e}t\'{e}s asymptotiques des groupes
lin\'{e}aires, GAFA, Geom. Funct. Anal. \textbf{7} (1997), 1-47.

\bibitem{Bourb}  Bourbaki, N. \textit{Groupes et alg\`{e}bres de Lie},
chapitres 7 et 8, Hermann, 1975.

\bibitem{BT}  Bruhat, F. ; Tits, J. Groupes r\'{e}ductifs sur un corps
local, I. Donn\'{e}es radicielles valu\'{e}es, \textit{Publ. math. IHES} 
\textbf{41} (1972).

\bibitem{D}  De la Harpe, P. \textit{Topics in Geometrical Group Theory},
Chicago University Press, (2000).

\bibitem{D}  De la Harpe, P., Free Groups in Linear Groups, \textit{%
L'Enseignement Math\'{e}matique}, \textbf{29} (1983), 129-144.

\bibitem{E}  Epstein, D.B.A., Almost all subgroups of a Lie group are free, 
\textit{J. of Algebra. }\textbf{19} (1971) 261-262.

\bibitem{Ghys}  Carri\`{e}re, Y., Ghys, E., Relations d'{\'e}quivalence moyennables sur les groupes de Lie, CRAS \textbf{300} (1985), no. 19, 677--680. 
 
\bibitem{kuranishi}  Kuranishi, M. On everywhere dense embedding of free
groups in Lie groups, \textit{Nagoya Math.} \textbf{J 2} (1951), 63-71.

\bibitem{PR}  Platonov, V., Rapinchuk, A., \textit{Algebraic Groups and
Number Theory}, Academic Press, 1994.

\bibitem{Rag}  Raghunathan, M.S. \textit{Discrete Subgroups of Lie Groups,
Ergebnisse der Mathematik und Ihrer Grenzgebiete}. Band \textbf{68} (1972).

\bibitem{T}  Tits, J. Free subgroups of Linear groups, \textit{Journal of
Algebra }\textbf{20}\textit{\ }(1972), 250-270.

\bibitem{Zim}  Zimmer, B.  Amenable actions and dense subgroups of Lie
groups, \textit{Journal of Functional analysis}, \textbf{72} (1987), no. 1, 58--64
 
\end{thebibliography}
\end{document}